\documentclass{svmult}

\smartqed
\usepackage{amsmath,amssymb}
\usepackage{ae}

\begin{document}

\title*{No Multiple Collisions for Mutually Repelling Brownian
Particles }
\author{Emmanuel C\'{e}pa and Dominique L\'{e}pingle}
\institute{{\small MAPMO, Universit\'{e} d'Orl\'{e}ans,
 \\ B.P.6759, 45067 Orl\'{e}ans Cedex 2, France} \\
 {\tt e-mail: Emmanuel.Cepa@univ-orleans.fr, dlepingl@univ-orleans.fr}}
\maketitle

\let\leq=\leqslant\let\le=\leqslant
\let\geq=\geqslant\let\ge=\geqslant

\def\N{\mathbb{N}}
\def\Z{\mathbb{Z}}
\def\Q{\mathbb{Q}}
\def\R{\mathbb{R}}
\def\C{\mathbb{C}}
\def\H{\mathbb{H}}
\def\E{\mathbb{E}}
\def\P{\mathbb{P}}
\def\L{\mathbb{L}}
\def\1{1\hspace{-1.2mm}\mbox{{\normalsize I}}}
\def\Dem{\medskip\noindent { Proof{}.}  }
\newcommand{\card}{\operatorname{card}}
\def\CQFD{\qed}


\newtheorem{Th}{Theorem}
\newtheorem{Def}[Th]{Definition}
\newtheorem{Lemme}[Th]{Lemma}
\newtheorem{Cor}[Th]{Corollary}
\newtheorem{Rmq}[Th]{Remark}
\newtheorem{Prop}[Th]{Proposition}
\newtheorem{Not}[Th]{Notation}
\newcommand\Section{\setcounter{equation}{0} \section}
\def\theequation{\thesection.\arabic{equation}}

 \noindent {\small{\bf Summary.} Although Brownian particles with 
 small mutual electrostatic 
 repulsion may collide, multiple collisions at
 positive time are always forbidden.}



\section{Introduction}
A three-dimensional Brownian motion $B_t= (B ^{1}_t , B ^{2}_t  , B ^{3}_t )$ does not hit the axis 
$\{x_1 = x_2=x_3 \}$ except possibly at time $0$. An easy proof is obtained by applying Ito's formula to 
$R_t \, = \, [ (B ^{1}_t  -B ^{2}_t )^2 +  (B ^{1}_t  -B ^{3}_t )^2 +  (B ^{2}_t  -B ^{3}_t )^2 ]$ and remarking 
that up to the multiplicative constant $3$ the process $R$ is the square of a two-dimensional 
Bessel process for which $\{0 \}$ is a polar state. This remark will be our guiding line in the sequel. 

We consider a filtered probability space $(\Omega , {\cal F} ,( {\cal F}
_t)_{t \geq 0} , \P)$ 
and for $N \geq 3$ the following system of stochastic differential equations
\begin{displaymath}  
dX_t^{i} \; = \;   dB_t^{i} \; + \; \displaystyle  
\lambda \sum_{1 \leq j \neq i \leq N} \frac 
{dt}{X_t^{i}-X_t^{j}}    \,
, \; i = 1 , 2 , \ldots , N 
\end{displaymath}
with boundary conditions
\begin{displaymath}  
X_t ^{1}\, \le \,  X_t ^{2}\, \le  \cdots \le \,  X_t ^{N} \, , \;    
\quad 0 \le  t < \infty  \, , 
\end{displaymath} 
and a random, ${\cal F} _0$-measurable, initial value satisfying
\begin{displaymath}  
X_0 ^{1}\, \le \,  X_0 ^{2}\, \le  \cdots \le \,  X_0 ^{N} \, . 
\end{displaymath}
Here $B_t= (B ^{1}_t , B ^{2}_t  ,  \ldots ,B ^{N}_t )$  denotes a standard 
$N$-dimensional $
( {\cal F} _t)$-Brownian motion and  $\lambda$ is a positive constant. This system has 
been extensively studied in 
\cite{Ch}, \cite{RS}, \cite{CL1}, \cite{BBCL}, \cite{CL2}, \cite{Fon}. 
For comments on the relationship between this system and the spectral analysis of Brownian 
matrices,
	and also conditioning of Brownian particles, we refer to the introduction and
the bibliography in \cite{CL2}.
	
When $\lambda \ge \displaystyle \frac{1}{2}$, establishing
strong existence and uniqueness is not difficult, because particles never
 collide, as proved in~\cite{RS}. The general case with arbitrary coupling strength 
 is investigated in \cite{CL1} and 
 it is proved in \cite{CL2} that collisions occur a.s. if and only if 
 $ 0< \lambda < \displaystyle \frac{1}{2}$. As for multiple collisions 
(three or more particles at the same location), it has been stated without proof in 
\cite{Sp} and \cite{CL3} that they are impossible. 
The proof we give below, with a Bessel process
unexpectedly coming in, is just an exercise on Ito's formula.

\section{A key Bessel process}

We consider for any $t \ge 0$
$$S_t \, = \,  \displaystyle  \sum _{j=1} ^N   \sum _{k=1} ^N   
(X ^{j}_t  -X ^{k}_t )^2 \; .$$

\begin{Th}
For any $\lambda > 0$, the process $S$ divided by the constant $2N$
is the square of a Bessel process with dimension 
$(N-1)(\lambda N +1)$.
\end{Th}

\Dem It is purely computational. Ito's formula provides for any $j \neq k$
$$ \begin{array}{lll}
(X ^{j}_t  -X ^{k}_t )^2 & = &  (X ^{j}_0  -X ^{k}_0 )^2 \, +  2 \displaystyle 
\int_0 ^t (X ^{j}_s  -X ^{k}_s ) d (B ^{j}_s  -B ^{k}_s ) \\

& &  +   2 \lambda \displaystyle  
 \sum_{1 \leq l \neq j \leq N}     \displaystyle \int_0^t      \frac {X_s^{j}-X_s^{k}} {X_s^{j}-X_s^{l}}   \,  ds 
\, +   2 \lambda \displaystyle  
 \sum_{1 \leq m \neq k \leq N}     \displaystyle \int_0^t      \frac {X_s^{k}-X_s^{j}} {X_s^{k}-X_s^{m}}  \,   ds \\ 
 && +  2\,t \, .
\end{array} $$
Adding the $N(N-1)$ equalities we get 

$$ \begin{array}{lll}
S_t & = & S_0 \, + \, 2 \displaystyle   \displaystyle  \sum_{j=1} ^N   
\sum_{k=1} ^N \int _0 ^t (X ^{j}_s  -X ^{k}_s ) d (B ^{j}_s  -B ^{k}_s ) \\
& &  + \,  4 \lambda \displaystyle  \sum_{j=1} ^N   \sum _{k=1} ^N 
 \sum _{1 \leq l \neq j \leq N}         \displaystyle \int _0 ^t          \frac {X_s^{j}-X_s^{k}} {X_s^{j}-X_s^{l}}    \, ds 
\,  + \, 2N(N-1)t \, .
\end{array} $$
But 
$$ \begin{array}{l}
 \displaystyle  \sum_{j=1} ^N   \sum_{k=1} ^N 
 \sum_{1 \leq l \neq j \leq N}         \displaystyle \int_0 ^t          
 \frac {X_s^{j}-X_s^{k}} {X_s^{j}-X_s^{l}}    \, ds \\
= \,   \displaystyle  \sum_{j=1} ^N   \sum_{k=1} ^N 

\displaystyle \sum_{1 \leq l \neq j \leq N}       \bigg [  \displaystyle \int _0 ^t          
\frac {X_s^{j}-X_s^{l}} {X_s^{j}-X_s^{l}}    \, ds

\,  + \, \displaystyle \int_0 ^t          \frac {X_s^{l}-X_s^{k}} {X_s^{j}-X_s^{l}}    \, ds  \bigg ] \\

= \; N^2 (N-1)t \, - \,  \displaystyle  \sum_{l=1} ^N   \sum_{k=1} ^N 
 \sum _{1 \leq j \neq l \leq N}         \displaystyle \int _0 ^t          \frac {X_s^{l}-X_s^{k}} {X_s^{l}-X_s^{j}}    \, ds \\

= \;       \displaystyle          \frac    {1}{2}     N^2 (N-1)t \, .

\end{array} $$

For the martingale term, we compute 

$$ \begin{array}{l}

\displaystyle  \sum _{j=1} ^N  ( \sum _{k=1} ^N   (X ^{j}_s  -X ^{k}_s) )^2 \\

= \;  \displaystyle  \sum _{j=1} ^N   \sum _{k=1} ^N   \sum _{l=1} ^N  (X ^{j}_s  -X ^{k}_s)  (X ^{j}_s  -X ^{l}_s)\\

= \;  \displaystyle  \sum _{j=1} ^N   \sum _{k=1} ^N   \sum _{l=1} ^N  (X ^{j}_s  -X ^{k}_s) ^2  \, + \, 

\displaystyle  \sum _{j=1} ^N   \sum _{k=1} ^N   \sum _{l=1} ^N (X ^{j}_s  -X ^{k}_s)  (X ^{k}_s  -X ^{l}_s)\\

= \;       \displaystyle          \frac    {N}{2}    S_s \, .

\end{array} $$

Let $B'$ be a linear Brownian motion independent of $B$. The process $C$ defined by :

$$C_t \, = \, \displaystyle \int _0 ^t  \1 _ {\{S_s > 0\}} \;   \displaystyle          \frac   {    \displaystyle  \sum _{j=1} ^N   \sum _{k=1} ^N     (X ^{j}_s  -X ^{k}_s) d B ^{j}_s } { \sqrt { \frac{N}{2} S_s} }
\; + \; \displaystyle \int _0 ^t  \1 _ {\{S_s = 0\}} dB'_s $$

is a linear Brownian motion and we have

$$ S_t \; = \; S_0  \, + \, 2 \displaystyle \int _0 ^t   \sqrt { 2 N S_s} dC_s 
 \,  + \, 2N(N-1)( \lambda N +1) t \, ,$$

which completes the proof.
\CQFD

\section{Multiple collisions are not allowed}
Since multiple collisions do not occur for Brownian particles without interaction, we can 
guess they do not
 either in case of mutual repulsion. Here is the proof.

\begin{Th}
For any $\lambda > 0$, multiple collisions cannot occur after time $0$.
\end{Th}
\Dem
i) For $3 \leq r \leq N$ and $1 \leq q \leq N-r+1$, let
\[
  \begin{array}{lll}
  I \, = \, \{ q, q+1, \ldots, q+r-1 \} \\ 
  S_t ^I\, = \,  \displaystyle  \sum _{j \in I}    
  \sum _{k \in I}   (X ^{j}_t  -X ^{k}_t )^2 \\   
  \tau ^I \, = \, \inf \{ t >0 : S_t ^I = 0 \}\;.
  \end{array}
\]   

ii) We first consider the initial condition $X_0$. From \cite{CL1}, Lemma 3.5, 
we know that for any $1 \le i < j \le N$ and any $t < \infty$, we have a.s.

$$  \displaystyle \int _0 ^t          \frac {du} {X_u^{j}-X_u^{i}}  \; < \; \infty \, . $$

Therefore for any $u >0$ there exists $0 < v < u$ such that 
$ X_v ^{1}\, < \,  X_v^{2}\, <  \cdots < \,  X_v ^{N} \; \mbox {a.s.} $
In order to prove $\P ( \tau ^I = \infty) =1$, we may thus assume $ X_0 ^{1}\, < \,  X_0^{2}\, <  \cdots < \,  X_0 ^{N} \; \mbox {a.s.} $, which implies for any $I$ that $S_0^I > 0$ and so $\tau ^I  > 0$ a.s. \\
iii) We know (\cite{RY}, XI, section 1) that $\{0\}$ is polar for the Bessel process $\sqrt { S_t} / \sqrt {2N}$, which means that  $\tau ^I = \infty$ a.s. for $I \, = \, \{ 1, 2, \ldots, N \}$. 
We will prove the same result for any $I$ by backward induction on $r = \card (I)$. 
Assume there are no $s$-multiple collisions for any $s>r$. Then

$$ \begin{array}{lll}
S_t ^I & = & S_0 ^I \, + \, 4 \displaystyle   \displaystyle  \sum _{j\in I}   \sum _{k \in I} 
\int _0 ^t (X ^{j}_s  -X ^{k}_s ) d B ^{j}_s   \\

& &  + \,  4 \lambda \displaystyle  \sum _{j \in I}    \sum _{k \in I} \sum _{l \notin I} 
         \displaystyle \int _0 ^t          \frac {X_s^{j}-X_s^{k}} {X_s^{j}-X_s^{l}}    \, ds 
\,  + \, 2r(r-1) ( \lambda r + 1) t \, .
\end{array} $$

We set for $n \in \N ^*$, $\tau ^I _n \, = \, \inf \{ t >0 : S_t ^I \le  \displaystyle \frac{1}{n} \}$.  For any $t \ge 0$,

$$ \begin{array}{lll}
\log S_{t \wedge \tau ^I _n} ^I & = & \log S_0 ^I \, + \, 4 \displaystyle   \displaystyle  \sum _{j\in I}   \sum _{k \in I} \int _0 ^ {t \wedge \tau ^I _n} 
\displaystyle \frac {X ^{j}_s  -X ^{k}_s }{S_s^I} d B ^{j}_s   \\

& &  + \,  2 \lambda \displaystyle  \sum _{j \in I}    \sum _{k \in I} \sum _{l \notin I} 
         \displaystyle \int _0 ^ {t \wedge \tau ^I _n}        \displaystyle \frac {(X ^{j}_s  -X ^{k}_s) }{S_s^I}   \bigg [ \frac {1}{X_s^{j}-X_s^{l}}   - \frac {1}{X_s^{k}-X_s^{l}}  \bigg ]  \, ds \\

&&
\,  + \, 2r[(r-1) ( \lambda r + 1) - 2]  \, \displaystyle \int _0 ^ {t \wedge \tau ^I _n}        
\displaystyle \frac {ds }{S_s^I}\\
 & > & - \infty \;.
\end{array} $$

>From the induction hypothesis we deduce that for $j,k \in I$ and $l  
\notin I$, a.s. on $\{ \tau ^I  < \infty \}$,  $ (X ^j  _ { \tau ^I } - X ^ l _ { \tau ^I }) (X ^k  _ { \tau ^I } - X ^ l _ { \tau ^I }) \, > \, 0$ and so

$$  \begin{array}{l}
\displaystyle \int _0 ^ {t \wedge \tau ^I }        
\displaystyle \frac {(X ^{j}_s  -X ^{k}_s) }{S_s}   
\bigg [ \frac {1}{X_s^{j}-X_s^{l}}   - \frac {1}{X_s^{k}-X_s^{l}}  \bigg ]  \, ds \\
= \; - \displaystyle \int _0 ^ {t \wedge \tau ^I }        
\displaystyle \frac {(X ^{j}_s  -X ^{k}_s)^ 2 }{S_s}   \frac {ds}{(X_s^{j}-X_s^{l} ) (X_s^{k}-X_s^{l})}   \\
\; - \infty \, .
\end{array} $$

The martingale $(M_n ,  {\cal F} _   {t \wedge \tau ^I _n} )_ {n \geq 1}$ 
defined by 
$$M_n \; = \;   4 \displaystyle   \displaystyle  \sum _{j\in I}   \sum _{k \in I} 
\int _0 ^ {t \wedge \tau ^I _n} 
\displaystyle \frac {X ^{j}_s  -X ^{k}_s }{S_s^I} d B ^{j}_s $$
has associated increasing process 
$A_n \, = \, 8 r     \displaystyle   \int _0 ^ {t \wedge \tau ^I _n}   
\displaystyle \frac {ds }{S_s^I} $. It follows that $M_n \, + \, \displaystyle \frac  {1}{4}
 [(r-1) ( \lambda r + 1) - 2] A_n$ 
either tends to a finite limit or to $+ \infty$ as $n$ tends to 
$+ \infty$. Then for any $t \ge 0$, $\log S_{t \wedge \tau ^I } ^I  \,  > \, - \infty$ and 
so $\P ( \tau ^I = \infty) =1$,
which completes the proof.
\CQFD

\section{ Brownian particles on the circle}
We now turn to the popular model of interacting Brownian 
particles on the circle (\cite{Sp}, \cite{CL2}). Consider the system of stochastic differential equations
\begin{displaymath}  
dX_t^{i} \; = \;   dB_t^{i} \; + \; \displaystyle  
\frac{\lambda}{2} \sum _{1 \leq j \neq i \leq N} 
\cot(\frac{X_t^i-X_t^j}{2}) 
dt    \,
, \; i = 1 , 2 , \ldots , N 
\end{displaymath}
with the boundary conditions
\begin{displaymath}  
X_t ^{1}\, \le \,  X_t ^{2}\, \le  \cdots \le \,  X_t ^{N} \le X_t^1\; +\; 2\pi \, , \;    
\quad 0 \le  t < \infty  \, .
\end{displaymath}
As expected we can prove there are no multiple collisions for the particles 
$Z_t^j\;=\; e^{i\,X_t^j} $ that live on the unit circle. The proof is more involved and 
will be deduced by approximation from the previous one.

\begin{Th}
Multiple collisions for the particles on the circle do not occur after time $0$ 
for any $\lambda > 0$.
\end{Th}

\noindent Sketch of the proof. For the sake of simplicity, we only deal with the 
$N$-collisions. Let
\[
 \begin{array}{lll}
 R_t & = &  \displaystyle  \sum _{j=1} ^N   \sum _{k=1} ^N   
\sin^2(\frac{X ^{j}_t  -X ^{k}_t}{2} ) \\
 \sigma_n & = & \inf\{t>0\,:\,R_t\leq \frac{1}{n}\} \;.
 \end{array}
\]
We apply Ito's formula to  $\log\,R_t$ and get
\[
  \log \,R_{t\wedge \sigma_n} \, = \, \log\,R_0 \,
    +\sum_{j=1}^N\,\int_0^{t\wedge \sigma_n}\,H_s^j\,dB_s^j \, + \, 
    \int_0^{t\wedge \sigma_n}L_s\,ds
    \]
for some continuous processes $H^{j}$ and $L$. We divide each integral into an integral
over $\{R_s\geq \frac{1}{2}\}$ and an integral over $\{R_s< \frac{1}{2}\}$. The first type 
integrals do not pose any problem. When $R_s<\frac{1}{2}$, we replace $X^j_s$ with 
\[
  Y^j_s\,=\, X^j_s \quad \mbox{ or } Y^j_s\,=\,X^j_s\,-\,2\pi
\]
  in such a way that for any $j,k$ we have $\mid Y_s^j-Y_s^k \mid<\pi/3$. The processes
  $H^{j}$ and $L$ have the same expressions in terms of $X$ or $Y$. 
  With this change of variables we may approximate $\sin x$ by $x$, 
   $\cos x$ by 1 and
  replace the trigonometric 
  functions by approximations of the linear ones which we have met in the previous sections.
We obtain that   
\[
  \log \,R_{t\wedge \sigma_n} \, = \, \log\,R_0
 \,+\,M_n \, + \, \displaystyle \frac  {1}{4}
 [(N-1) ( \lambda N + 1) - 2] A_n \,+
  \,\int_0^{t\wedge \sigma_n}\,D_s\,ds
  \]
where $M_n$ is a martingale with associated increasing process $A_n$ and $D$ 
is a.s. a locally integrable  process. Details are left to the reader as well as the case of an 
 arbitrary subset $I$  
 like those in Section 3. 
\CQFD

\end{document}